\documentclass[a4paper]{article}
\usepackage{amsfonts}
\usepackage{amssymb}
\usepackage{latexsym}
\newtheorem{theorem}{Theorem}[section]

\newtheorem{remark}[theorem]{Remark}
\newtheorem{example}[theorem]{Example}

\begin{document}

\parbox{1mm}

\begin{center}
{\bf {\sc \Large Some special functions identities arising from
commuting operators}}
\end{center}

\vskip 12pt

\begin{center}
{\bf M\'aria Hutn\'ikov\'a and Ondrej Hutn\'{\i}k}
\end{center}

\vskip 24pt

\hspace{5mm}\parbox[t]{11cm}{\fontsize{9pt}{0.1in}\selectfont\noindent{\bf
Abstract.} Commuting is an important property in many cases of
investigation of properties of operators as well as in various
applications, especially in quantum physics. Using the observation
that the generalized weighted differential operator of order $k$
and the weighted Hardy-type operator commute we derive a number of
new and interesting identities involving some functions of
mathematical physics.}

\vskip 24pt

\section{Introduction}

When studying two operators $P,Q$ of quantum theory it is crucial
whether the relation $$[P,Q]=0, \quad \textrm{where }\,\,
[P,Q]=(-\imath\hbar)^{-1}(PQ-QP),$$ is fulfilled or not. Such
physical quantities, for which there are uniquely assigned
operators $P,Q$, are simultaneously measurable if and only if
$[P,Q]=0$. Equivalently, we say that the \textit{operators $P,Q$
commute} if the commutator $[P,Q]=0$, i.e. $PQ=QP$. There are many
known operators for which the commutation relation is fulfilled,
e.g. the class of normal operators. On the other side,
non-commutative operators are source of some interesting stories
in mathematics and physics, e.g. the famous Heisenberg uncertainty
principle saying that it is impossible to know the momentum and
position of a particle simultaneously (see~\cite{neumann}).

In this paper we are interested in two operators, namely the
\textit{Hardy-type operator} $\mathbf{H}_{w}$ defined by
$$\mathbf{H}_{w}f(x)=\frac{1}{w_1(x)}\int_{\alpha}^{x}w(t)f(t)\,\textrm{d}t, \quad x>0,$$
where $-\infty \leq \alpha\leq \infty$,
$w_1(x)=\int_{\alpha}^{x}w(t)\,\textrm{d}t$ and $w,f$ are real
measurable, locally integrable functions, and the
\textit{generalized weighted differential operator}
$\mathbf{D}_{k}^{w}$ of order $k$ given by
$$\mathbf{D}_{k}^{w} f(x) = \left(\frac{w_1(x)}{w(x)} \frac{d}{dx}\right)^{k} f(x) =
\left(\frac{w_1(x)}{w(x)}
\frac{d}{dx}\right)^{k-1}\left(\frac{w_1(x)}{w(x)} f'(x)\right),$$
$k=1, 2,\dots$. Note that operators $\mathbf{H}_w$ are of the
interest in various functional spaces mainly in connection with
Hardy's inequality, cf.~\cite{persson}. The case $\alpha=0$ and
$w\equiv \textrm{const}$ corresponds to the Hardy's averaging
operator (or, Hardy's arithmetic mean operator), therefore we
propose to say $\mathbf{H}_w$ the Hardy-type operator. On the
other hand, operators $\mathbf{D}_k^w$ appear in the theory
of~differential equations (see e.g.~\cite{mojsej}
and~\cite{ohriska}). If $w_1(x)/w(x)=1$, the resulting $k$-th
order differential operator will be simply denoted by
$$\mathbf{D}_{k}f(x)=\frac{d^k}{dx^k} f(x), \quad k=1, 2,\dots,$$
as usual. A generalization of weighted Hardy's averaging operator
is provided in~\cite{hutnik}, where the general mean-type
inequality involving such operator is investigated.

The first result of this paper is the observation that the
operators $\mathbf{D}_{k}^{w}$ and $\mathbf{H}_{w}$ commute. Using
this fact in Section~3 we establish few interesting new identities
involving some special functions of mathematical physics.

\section{A note on commutation relation}

Now we describe an easy observation on commutativity of operators
$\mathbf{D}_k^w$ and $\mathbf{H}_w$. 
For the sake of brevity let us replace $\mathbf{H}_{w}h(x)$ by
$H(x)$. Then we have
\begin{equation}\label{eq2thm1}
H(x)=\frac{1}{w_1(x)}\int_{\alpha}^{x} w(t)h(t)\,\textrm{d}t.
\end{equation} Differentiating the equality~(\ref{eq2thm1}) 
and then using integration by parts with $w_{1}(\alpha)=0$ we
obtain the identity
$$H'(x)=\frac{w(x)\int_{\alpha}^{x}w_{1}(t)h'(t)\,\textrm{d}t}{w_{1}(x)^{2}},$$
which may be written as follows
\begin{eqnarray}\label{eq3thm1}
\frac{w_{1}(x)}{w(x)} H'(x) & = & \frac{1}{w_{1}(x)}
\int_{\alpha}^{x} w_{1}(t) h'(t)\,\textrm{d}t \nonumber \\ & = &
\frac{1}{w_{1}(x)}\int_{\alpha}^{x} w(t) \left(
{\frac{w_{1}(t)}{w(t)}} h'(t)\right)\,\textrm{d}t,
\end{eqnarray}
or in the terms of weighted differential operator as
$$\mathbf{D}_{1}^{w}H(x)=\frac{1}{w_1(x)}\int_{\alpha}^{x}w(t)[\mathbf{D}_{1}^{w}h](t)\,\textrm{d}t.$$
Induction on $k$ gives
\begin{equation}\label{eqthm1}
\mathbf{D}_{k}^{w}\mathbf{H}_{w}h(x) =
\frac{1}{w_1(x)}\int_{\alpha}^{x}w(t)
\mathbf{D}_{k}^{w}h(t)\,\textrm{d}t, \qquad k=1, 2, \ldots \ .
\end{equation} Summarizing the above we get

\begin{theorem}\label{theorem}
For each $k\in\mathbb{N}$ operators $\mathbf{D}_{k}^{w}$ and
$\mathbf{H}_{w}$ commute, i.e.
\begin{equation}\label{eqthm1'}
\mathbf{D}_{k}^{w}\mathbf{H}_{w}=
\mathbf{H}_{w}\mathbf{D}_{k}^{w}.
\end{equation}
\end{theorem}

\begin{remark}\rm
The observation that operators $\mathbf{D}_{k}^{w}$ and
$\mathbf{H}_{w}$ commute for each $k\in\mathbb{N}$ may be also
verified using the well-known fact that the operators
$$\mathbf{A}f(x) = f(x) + \frac{w_1(x)}{w(x)} f'(x)$$ and
$\mathbf{H}_w$ are the inverses of each other.
\end{remark}

\begin{remark}\rm\label{remCL}
Note that for the weighted Hardy's averaging operator (the case
$\alpha=0$) it is enough to consider the operators
$$\mathbf{C}f(x) = e^{-x} \int_{-\infty}^{x} e^t
f(t)\,\textrm{d}t,$$ and $\mathbf{L}f(x) = f\left(\ln
w_1(x)\right)$ for positive function $w$. Taking into account the
relations
$$\mathbf{L}\mathbf{C} = \mathbf{H}_w\mathbf{L}, \quad
\mathbf{L}\mathbf{D}_k = \mathbf{D}_k^w \mathbf{L},$$ for each
$k\in \mathbb{N}$, we have the commutativity of operators
$\mathbf{C}$ and $\mathbf{D}_k$, and hence $\mathbf{H}_w$ and
$\mathbf{D}_k^w$.
\end{remark}

Further, let us denote $w_{k}(x) = \int_{\alpha}^{x}
w_{k-1}(x)\,\textrm{d}x$, where $w_0(x)=w(x)$, and $$r_{k}(x) =
\frac{w_{k}^{2}(x)}{w_{k-1}(x) w_{k+1}(x)},\quad k =1,2, \ldots.$$
Using $w_1(\alpha)=0$ and the relation
$$w_1^2(x)\,\frac{d}{dx}\left(\frac{1}{w_1(x)}\int_\alpha^x w(t)f(t)\,\textrm{d}t\right) = w(x)\int_\alpha^x w_1(t)f'(t)\,\textrm{d}t,$$
we get
$$r_{1}(x) {\frac{d}{dx}} \, \mathbf{H}_{w}h(x)=
\frac{1}{w_{2}(x)}\int_{\alpha}^{x}
w_{1}(t)\mathbf{D}_{1}h(t)\,\textrm{d}t.$$ When we repeat this
argument $k$-times with $w_{\ell + 1}(x)= \int_{0}^{x}
w_{\ell}(t)\,\textrm{d}t$, for $\ell = 1, 2, \ldots , k$, we
obtain {\setlength\arraycolsep{2pt}
\begin{eqnarray}\label{eq0thm1}
& {\phantom{1}}& r_{k}(x) {\frac{d}{dx}} \left(r_{k-1}(x)
{\frac{d}{dx}}\left( \ldots r_{1}(x) {\frac{d}{dx}} \,
\mathbf{H}_{w}h(x)\right)\right) \nonumber \\ & = &
\frac{1}{w_{k+1}(x)}\int_{\alpha}^{x}
w_{k}(t)\mathbf{D}_{k}h(t)\,\textrm{d}t, \quad k=1,2, \ldots \ ,
\end{eqnarray}}which yields the following

\begin{theorem}\label{theorem2}
In the sense of the above notation
\begin{equation}\label{eq0thm1'}
\mathbf{R}_{k}\mathbf{H}_{w} =\mathbf{H}_{w_{k}}\mathbf{D}_{k},
\end{equation} where the operator $\mathbf{R}_{k}$ is given by
$$\mathbf{R}_{k}\nu(x)=r_{k}(x) \,
\frac{d}{dx} \left( r_{k-1}(x) \, \frac{d}{dx}\left( \dots
r_{1}(x) \, \frac{d}{dx} \, \nu(x)\right)\right).$$
\end{theorem}

\begin{remark}\rm\label{remRk}
The operator $\mathbf{R}_k$ will be called the
\textit{quasi-differential operator of order $k$} because it
appears in connection with differential equations with
quasi-derivatives, see e.g.~\cite{mojsej}. Clearly, if $r_k(x)=1$
for each $k\in\mathbb{N}$, then $\mathbf{R}_k$ reduces to
$\mathbf{D}_k$ and it is easy to observe that the usual $k$-th
differential operator $\mathbf{D}_k$ commutes with $\mathbf{H}_w$
for each $k\in\mathbb{N}$.
\end{remark}

Note that the integral identities~(\ref{eqthm1'})
and~(\ref{eq0thm1}) may be equivalently written as differential
identities in the form
\begin{equation}\label{dif1}
\mathbf{D}_{1}\left(w_1(x)\mathbf{D}_{k}^{w}\mathbf{H}_{w}h(x)\right)
= w(x)\mathbf{D}_{k}^{w}\left(\frac{1}{w(x)}
\mathbf{D}_{1}\Bigl(w_1(x)\mathbf{H}_{w}h(x)\Bigr)\right),
\end{equation} and
\begin{equation}\label{dif2}
\mathbf{D}_{1}\left(w_{k+1}(x)
\mathbf{R}_{k}\mathbf{H}_{w}h(x)\right) = w_{k}(x)\mathbf{D}_{k}
\left(\frac{1}{w(x)}\mathbf{D}_{1}\Bigl(w_{1}(x)\mathbf{H}_{w}h(x)\Bigr)\right),
\end{equation} respectively. Indeed, the term
$$\mathbf{H}_{w}h(x)=\frac{1}{w_1(x)}\int_{\alpha}^{x}
w(t)h(t)\,\textrm{d}t$$ is equivalent to
$$\frac{d}{dx}\Bigl(w_{1}(x)\mathbf{H}_{w}h(x)\Bigr)=w(x)h(x).$$
Now the identity~(\ref{eqthm1'}) may be rewritten as
\begin{equation}\label{eqdif}
\frac{d}{dx}\Bigl(w_{1}(x)\mathbf{D}_{k}^{w}\mathbf{H}_{w}h(x)\Bigr)=
w(x)\mathbf{D}_{k}^{w}h(x).
\end{equation}
If we substitute
\begin{equation}\label{eqdif3}
h(x)=w^{-1}(x)\frac{d}{dx}\Bigl(w_{1}(x)[\mathbf{H}_{w}h](x)\Bigr)
\end{equation}
into~(\ref{eqdif}), we obtain the following differential operator
identity
$$\frac{d}{dx}\left(w_1(x)\mathbf{D}_{k}^{w}\mathbf{H}_{w}h(x)\right)
= w(x)\mathbf{D}_{k}^{w}\left(\frac{1}{w(x)}
\frac{d}{dx}\Bigl(w_1(x)\mathbf{H}_{w}h(x)\Bigr)\right),$$ which
corresponds to~(\ref{dif1}). Similarly, the integral
identity~(\ref{eq0thm1}) is equivalent to the following
differential identity {\setlength\arraycolsep{2pt}
\begin{eqnarray*}
& {\phantom{1}}& \frac{d}{dx}\left(w_{k+1}(x)\,r_{k}(x) \, {\frac{
d }{dx}} \left( r_{k-1}(x) \, {\frac{ d }{dx}}\left( \dots
r_{1}(x) \, {\frac{ d }{dx}} \,
\mathbf{H}_{w}h(x)\right)\right)\right) \\ & = &
w_{k}(x)\mathbf{D}_{k}h(x),
\end{eqnarray*}}which may be rewritten
by the use of substitution~(\ref{eqdif3}) and the
quasi-differential operator $\mathbf{R}_{k}$ to desired
form~(\ref{dif2}).

\section{Applications to special functions}

The obtained relations~(\ref{eqthm1'}) and~(\ref{eq0thm1'})
involve functions defined as integrals of upper limit. Such
functions usually appear in the theory of special functions.
Special functions, also denoted as functions of mathematical
physics, have important applications in many areas of mathematics,
science and engineering. Using the observation on commuting
operators we now present few applications of the above obtained
relations to some special functions to derive some interesting
identities and representations for them. As far as we know these
identities are not known in the available literature, e.g.
in~\cite{abramovitz} and~\cite{rainville}.

In the following few examples we consider $\alpha = 0$ (the
Hardy's averaging operator) and polynomial weights.

\begin{example}\rm
Considering the functions $h(t)=e^{-t}$ and $w(t)=t^{a-1}$ for
$a>0$ we get $w_{1}(x)/w(x)=x/a$. Then
$$\mathbf{H}_{w}h(x)=\frac{1}{\int_{0}^{x}t^{a-1}\,\textrm{d}t}\int_{0}^{x}t^{a-1}e^{-t}\,\textrm{d}t
=\frac{a}{x^{a}} \gamma(a,x),$$ where $\gamma(a,x)$ denotes the
(lower) incomplete gamma function (see~\cite{abramovitz}). For the
generalized $k$-th derivative of $\mathbf{H}_{w}h(x)$ we get
$$\mathbf{D}_{k}^{w}\mathbf{H}_{w}h(x)=
\frac{1}{a^{k-1}}\left(x\frac{d}{dx}\right)^{k}\Bigr(x^{-a}\gamma(a,x)\Bigr),
\qquad k=1,2,\dots \ .$$ On the other hand, the direct calculation
of $\mathbf{H}_{w}\mathbf{D}_{k}^{w}h(x)$ yields
{\setlength\arraycolsep{2pt}
\begin{eqnarray*}
\mathbf{H}_{w}\mathbf{D}_{k}^{w}h(x) & = &
\frac{1}{x^{a}a^{k-1}}\int_{0}^{x}t^{a-1}\left(t\frac{d}{dt}\right)^{k}\left(e^{-t}\right)\,\textrm{d}t
\\ & = & \frac{1}{x^{a}a^{k-1}}\sum_{i=1}^{k}(-1)^{k-i+1}
S_2(k,k-i+1)\int_{0}^{x}t^{a+k-i}e^{-t}\,\textrm{d}t \\
& = &
\frac{1}{x^{a}a^{k-1}}\sum_{i=1}^{k}(-1)^{k-i+1}S_2(k,k-i+1)\gamma(a+k-i+1,x),
\end{eqnarray*}}where $S_2(n,m)$ is the Stirling number of the second kind
(see~\cite{abramovitz}). By Theorem~\ref{theorem} we obtain the
relation
$$\left(x\frac{d}{dx}\right)^{k}\Bigr(x^{-a}\gamma(a,x)\Bigr)=
\frac{1}{x^{a}}\sum_{i=1}^{k}(-1)^{k-i+1}S_2(k,k-i+1)\gamma(a+k-i+1,x).$$
For $k=1$ we get the well known relation $\gamma(a+1, x) =
a\gamma(a,x) - x^{a}e^{-x}$. Similarly, for $k=2$ we have
$\gamma(a+2,x) = a(a+1)\gamma(a,x)-(a+1+x)x^{a}e^{-x}$.
\end{example}

\begin{remark}\rm
In the literature (see~\cite{abramovitz}) there exists the
following representation
\begin{equation}\label{repgamma}
\gamma(a+n,x)=(-1)^{n}x^{a+n} \frac{d^n}{dx^n}(x^{-a}\gamma(a,x)).
\end{equation} Therefore,
$$\left(x\frac{d}{dx}\right)^{k}\Bigr(x^{-a}\gamma(a,x)\Bigr)=
\sum_{i=1}^{k}S_2(k,k-i+1)x^{k-i+1}\frac{d^{k-i+1}}{dx^{k-i+1}}(x^{-a}\gamma(a,x)).$$
\end{remark}

\begin{example}\rm
Let us consider $w(t)=t^n$, $n\in \mathbb{N}_0$, and
$h(t)=\psi_{0}(t)$, where $\psi_{0}(t)$ is the digamma function.
Using the same method as in the previous example, we get
$$\mathbf{D}_{k}^{w}\mathbf{H}_{w}h(x)=
\frac{1}{(n+1)^{k-1}} \left(x \frac{d}{dx}\right)^{k}
\left(\frac{1}{x^{n+1}} \int_{0}^{x}t^{n}
\psi_{0}(t)\,\textrm{d}t\right),$$ and
$$\mathbf{H}_{w}\mathbf{D}_{k}^{w}h(x)=
\frac{1}{x^{n+1}(n+1)^{k-1}} \int_{0}^{x}t^{n} \left(t
\frac{d}{dt}\right)^{k} \psi_{0}(t)\,\textrm{d}t.$$ Since $\left(t
\frac{d}{dt}\right)^{k} \psi_{0}(t) =
\sum\limits_{i=1}^{k}S_2(k,i)\,t^i\psi_i(t)$, where
$\psi_n(t)=\frac{d^n}{dt^n}\psi_0(t)$ is the polygamma function,
then $$\mathbf{H}_{w}\mathbf{D}_{k}^{w}h(x)=
\frac{1}{x^{n+1}(n+1)^{k-1}}
\sum_{i=1}^{k}S_2(k,i)\int_{0}^{x}t^{n+i}\psi_{i}(t)\,\textrm{d}t.$$
Thus,
$$\left(x \frac{d}{dx}\right)^{k}
\left(\frac{1}{x^{n+1}} \int_{0}^{x}t^{n}
\psi_{0}(t)\,\textrm{d}t\right)=\frac{1}{x^{n+1}}
\sum_{i=1}^{k}S_2(k,i)\int_{0}^{x}t^{n+i}\psi_{i}(t)\,\textrm{d}t.$$
Finally, for $k=1$ we get {\setlength\arraycolsep{2pt}
\begin{eqnarray*}
\int_{0}^{x}t^{n+1} \psi_{1}(t)\,\textrm{d}t & = & x^{n+2}
\frac{d}{dx}\left(\frac{1}{x^{n+1}} \int_{0}^{x}t^{n}
\psi_{0}(t)\,\textrm{d}t\right)
\\ & = &
x^{n+1}\psi_{0}(x)-(n+1)\int_{0}^{x}t^{n}
\psi_{0}(t)\,\textrm{d}t,
\end{eqnarray*}}which gives a formula for the $(n+1)$-th moment of the
trigamma function by the use of digamma function.
\end{example}

\begin{example}\rm
Let $w(t)=t^m$, $m\in \mathbb{N}_0$, $h(t)=B_n(t)$, $n\in
\mathbb{N}_0$, where $B_n(t)$ is the $n$-th Bell polynomial
(see~\cite{abramovitz}). Then
$$\mathbf{D}_{k}^{w}\mathbf{H}_{w}h(x)= \frac{1}{(m+1)^{k}}
\left(x \frac{d}{dx}\right)^{k} \left(\frac{1}{x^{m+1}}
\int_{0}^{x}t^{m} B_{n}(t)\,\textrm{d}t\right),$$ and
$$\mathbf{H}_{w}\mathbf{D}_{k}^{w}h(x)=
\frac{1}{x^{m+1}(m+1)^{k}} \int_{0}^{x}t^{m} \left(t
\frac{d}{dt}\right)^{k} B_{n}(t)\,\textrm{d}t.$$ Using the
explicit formula $$B_n(t)=\sum_{i=0}^{n} S_2(n,i)\,t^i$$ we get
$$\frac{1}{x^{m+1}(m+1)^{k}} \int_{0}^{x}t^{m} \left(t
\frac{d}{dt}\right)^{k} B_{n}(t)\,\textrm{d}t =
\frac{1}{(m+1)^{k}} \sum_{i=1}^{n}
\frac{i^k}{m+i+1}\,S_2(n,i)\,x^i.$$ Thus,
$$\left(x \frac{d}{dx}\right)^{k} \left(\frac{1}{x^{m+1}}
\int_{0}^{x}t^{m} B_{n}(t)\,\textrm{d}t\right) = \sum_{i=1}^{n}
\frac{i^k}{m+i+1}\,S_2(n,i)\,x^i.$$ Specially, when $m=0$ and
$k=1$ we get the relation for integral of the $n$-th Bell
polynomial
$$\int_{0}^{x}B_n(t)\,\textrm{d}t = xB_{n}(x)-\sum_{i=1}^{n}
\frac{i}{i+1}\,S_2(n,i)\,x^i.$$
\end{example}

\begin{example}\rm
Consider $w(t)=1$ and $h(t)=(t^2-1)^k$, $k\in \mathbb{N}$. Then we
have
$$w_{k}(x)=\frac{x^{k}}{k!}, \quad \textrm{and} \quad
r_{k}(x)=\frac{w^{2}_{k}(x)}{w_{k-1}(x)w_{k+1}(x)}=\frac{k+1}{k},
\quad k=1,2,\dots .
$$ A direct calculation yields
$$\mathbf{R}_{k}\mathbf{H}_{w}h(x)=(k + 1)
\frac{d^k}{dx^k} \left(\frac{1}{x}
\int_{0}^{x}(t^2-1)^k\,\textrm{d}t\right),$$ and
$$\mathbf{H}_{w_{k}}\mathbf{D}_{k}h(x)=\frac{k+1}{x^{k+1}}
\int_{0}^{x} t^k \frac{d^k}{dt^k}(t^2-1)^k\,\textrm{d}t =
\frac{2^{k}(k+1)!}{x^{k+1}}\int_{0}^{x} t^k P_k(t)\,\textrm{d}t,$$
where the Rodrigues representation $$P_k(x) = \frac{1}{2^k k!}
\frac{d^k}{dx^k}(x^2-1)^k$$ of Legendre polynomials $P_k(x)$ has
been used (see~\cite{abramovitz}). Then by Theorem~\ref{theorem2}
we get
$$\int_{0}^{x} t^k P_k(t)\,\textrm{d}t=\frac{x^{k+1}}{2^k k!}\frac{d^k}{dx^k} \left(\frac{1}{x}
\int_{0}^{x}(t^2-1)^k\,\textrm{d}t\right).$$ However, direct
calculations yield
$$\frac{1}{x} \int_{0}^{x}(t^2-1)^k\,\textrm{d}t = \sum_{i=0}^{k}\frac{(-1)^i{k \choose
i}}{2(k-i)+1}\,x^{2(k-i)},$$ and therefore we have
$$\frac{d^k}{dx^k} \left(\frac{1}{x}
\int_{0}^{x}(t^2-1)^k\,\textrm{d}t\right) = \sum_{i=0}^{k}
\frac{(-1)^i{k\choose i}(2(k-i))!}{(2(k-i)+1)(k-2i)!}\,x^{k-2i},$$
from which follows the identity $$\int_{0}^{x} t^k
P_k(t)\,\textrm{d}t = \frac{1}{2^k k!} \sum_{i=0}^{k}
\frac{(-1)^i{k\choose i}(2(k-i))!}{(2(k-i)+1)(k-2i)!}\,x^{k-2i}.$$
\end{example}

\begin{example}\rm
As in the previous example, put $w(t)=1$ and
$h(t)=\frac{2}{\sqrt{\pi}} e^{-t^2}$. Then
$$\mathbf{R}_{k}\mathbf{H}_{w}h(x)=(k+1)
\frac{d^k}{dx^k}\left(\frac{2}{\sqrt{\pi}x}
\int_{0}^{x}e^{-t^2}\,\textrm{d}t\right)=(k+1)
\frac{d^k}{dx^k}\left(\frac{\textrm{erf}(x)}{x}\right),$$ where
$\textrm{erf}(x)$ is the Gauss error function. On the other hand
{\setlength\arraycolsep{2pt}
\begin{eqnarray*}
\mathbf{H}_{w_{k}}\mathbf{D}_{k}h(x) & = & \frac{k+1}{x^{k+1}}
\int_{0}^{x} t^k \frac{d^k}{dt^k}\left(\frac{2}{\sqrt{\pi}}
e^{-t^2}\right)\,\textrm{d}t \\ & = &
\frac{2(-1)^{k}}{\sqrt{\pi}}\frac{k+1}{x^{k+1}} \int_{0}^{x} t^k
H_{k}(t) e^{-t^2}\,\textrm{d}t,
\end{eqnarray*}}where $H_k$ is a Hermite
polynomial (see~\cite{abramovitz}). Therefore we get the relation
$$\frac{2}{\sqrt{\pi}}\int_{0}^{x} t^k H_{k}(t) e^{-t^2}\,\textrm{d}t=(-1)^{k}x^{k+1}
\frac{d^k}{dx^k}\left(\frac{\textrm{erf}(x)}{x}\right).$$
Moreover, using the MacLaurin series representation
$$\textrm{erf}(x)=\frac{2}{\sqrt{\pi}}\sum_{n=0}^{\infty} \frac{(-1)^n}{n!(2n+1)}\,x^{2n+1},$$
we get $$\int_{0}^{x} t^k H_{k}(t)
e^{-t^2}\,\textrm{d}t=\sum_{n=0}^{\infty}\frac{(-1)^{n+k}(2n)!}{(2n+1)n!(2n-k)!}\,x^{2n-k},
\quad x>0.$$
\end{example}

\begin{example}\rm
Let $h(t)=t^ke^{-t}$, $k\in\mathbb{N}_0$ and $w(t)=1$. Then
$$\mathbf{R}_{k}\mathbf{H}_{w}h(x)=(k+1)
\frac{d^k}{dx^k}\left(\frac{1}{x} \int_{0}^{x}t^k
e^{-t}\,\textrm{d}t\right)= (k+1)
\frac{d^k}{dx^k}\left(\frac{\gamma(k+1,x)}{x}\right),$$ and
$$\mathbf{H}_{w_{k}}\mathbf{D}_{k}h(x) = \frac{k+1}{x^{k+1}}
\int_{0}^{x} t^k \frac{d^k}{dt^k}\left(t^k
e^{-t}\right)\,\textrm{d}t = \frac{(k+1)!}{x^{k+1}} \int_{0}^{x}
t^k e^{-t}L_k(t)\,\textrm{d}t,$$ where the Rodrigues
representation
$$L_n(x)=\frac{e^x}{n!}\frac{d^n}{dx^n}(x^n e^{-x})$$ for the Laguerre
polynomials $L_k$ is used. Then $$\int_{0}^{x} t^k
e^{-t}L_k(t)\,\textrm{d}t = \frac{x^{k+1}}{k!}
\frac{d^k}{dx^k}\left(\frac{\gamma(k+1,x)}{x}\right).$$ According
to~(\ref{repgamma}) we have $$\frac{\gamma(k+1,x)}{x} = (-1)^k
x^{k}\frac{d^k}{dx^k}\left(\frac{\gamma(1,x)}{x}\right) = (-1)^k
x^{k}\frac{d^k}{dx^k}\left(\frac{1-e^{-x}}{x}\right),$$ and
therefore $$\frac{x^{k+1}}{k!}
\frac{d^k}{dx^k}\left(\frac{\gamma(k+1,x)}{x}\right) = (-1)^k
x^{k+1}\sum_{i=0}^{k}{k \choose
i}\frac{x^i}{i!}\frac{d^{k+i}}{dx^{k+i}}\left(\frac{1-e^{-x}}{x}\right).$$
Since
$$\frac{d^{k+i}}{dx^{k+i}}\left(\frac{1-e^{-x}}{x}\right)=\frac{(-1)^{k+i}(k+i)!}{x^{k+i+1}}
+ (-1)^{k+i+1}\sum_{j=0}^{k+i} T(k+i,j) \frac{e^{-x}}{x^{j+1}},$$
where $T(n,k)$ are the permutation coefficients giving number of
permutations of $n$ things $k$ at a time, then
$$\int_{0}^{x} t^k e^{-t}L_k(t)\,\textrm{d}t = \sum_{i=0}^{k}\frac{(-1)^{i+1}}{i!}{k \choose i}
\left[e^{-x}\sum_{j=0}^{k+i}T(k+i,j)x^{k+i-j} - (k+i)!\right].$$
\end{example}

In the next three examples we derive special function relations
involving $\alpha=+\infty$, resp. $\alpha=-\infty$, and
non-polynomial weights.

\begin{example}\rm
Let $\alpha = +\infty$ and $w(t)=e^{-t}$. Then $w_k(t)=e^{-t}$ for
each $k\in\mathbb{N}$, and therefore $r_k(t)=1$ for each
$k\in\mathbb{N}$. Choosing $h(t)=t^{a-1}$, $a>0$, according to
Remark~\ref{remRk} we have {\setlength\arraycolsep{2pt}
\begin{eqnarray*}
e^x \int_{x}^{\infty} e^{-t}\,\frac{d^k}{dt^k}
(t^{a-1})\,\textrm{d}t & = & \frac{d^k}{dx^k} \left(e^x
\int_{x}^{\infty} e^{-t}t^{a-1}\,\textrm{d}t\right) \\ & = &
\frac{d^k}{dx^k} (e^x \Gamma(a,x)),
\end{eqnarray*}}where $\Gamma(a,x) = \Gamma(a)-\gamma(a,x) = \int_{x}^{\infty}
e^{-t}t^{a-1}\,\textrm{d}t$ is the (upper) incomplete gamma
function (see~\cite{abramovitz}). Since
$$\frac{d^k}{dx^k} (e^x \Gamma(a,x)) = e^x \sum_{i=0}^{k}\frac{d^i}{dx^i}\Gamma(a,x),$$
then we finally get the relation
$$\sum_{i=0}^{k}\frac{d^i}{dx^i}\Gamma(a,x) = (a-1)_k\, \Gamma(a-k,x), \quad
a>k,$$ where $(x)_n=x(x-1)\dots(x-(n-1))$ is the falling
factorial. Observe that the case $k=1$ in the last relation
corresponds to the well-known relation
$$\Gamma(a+1,x)=a\Gamma(a,x)+x^ae^{-x}.$$
\end{example}

\begin{example}\rm
Let $\alpha=-\infty$. Consider functions $w(t)=e^t$ and
$f(t)=\exp\left(-\frac{(t-m)^2}{2\sigma^2}-t\right)$, where
$-\infty < m < \infty$, and $0<\sigma<\infty$. Then by
Remark~\ref{remRk} we have
$$\frac{d^k}{dx^k} \left(e^{-x} \int_{-\infty}^{x}
e^{-\frac{(t-m)^2}{2\sigma^2}}\textrm{d}t\right) = e^{-x}
\int_{-\infty}^x e^{t}\,\frac{d^k}{dt^k}
e^{-\left(\frac{(t-m)^2}{2\sigma^2}+t\right)}\,\textrm{d}t.$$ Note
that the integral $\int_{-\infty}^{x}
e^{-\frac{(t-m)^2}{2\sigma^2}}\textrm{d}t$ is the normal
distribution function with mean $m$ and standard deviation
$\sigma$. Since
$$\int_{-\infty}^{x} e^{-\frac{(t-m)^2}{2\sigma^2}}\textrm{d}t =
\sigma
\sqrt{\frac{\pi}{2}}\left(1+\textrm{erf}\left(\frac{x-m}{\sqrt{2}\sigma}\right)\right),$$
see~\cite{abramovitz}, then by Leibnitz rule we have
$$\frac{d^k}{dx^k} \left(e^{-x} \int_{-\infty}^{x}
e^{-\frac{(t-m)^2}{2\sigma^2}}\textrm{d}t\right) = \sigma
\sqrt{\frac{\pi}{2}}e^{-x}\sum_{i=0}^{k}(-1)^i{k \choose
i}\frac{d^{k-i}}{dx^{k-i}}\textrm{erf}\left(\frac{x-m}{\sqrt{2}\sigma}\right).$$
Using the relation
$$\frac{d^{n+1}}{dx^{n+1}}\textrm{erf}(x) = (-1)^n\frac{2}{\sqrt{\pi}}H_n(x)e^{-x^2}, \,\, n\in\mathbb{N}_0,$$
(see~~\cite{abramovitz}) where $H_n$ is the Hermite polynomial,
yields {\setlength\arraycolsep{2pt}
\begin{eqnarray*}
& {\phantom{1}}& \frac{d^k}{dx^k} \left(e^{-x} \int_{-\infty}^{x}
e^{-\frac{(t-m)^2}{2\sigma^2}}\textrm{d}t\right) \\ & = &
(-1)^k\sigma
\sqrt{\frac{\pi}{2}}e^{-x}\left[\textrm{erf}\left(\frac{x-m}{\sqrt{2}\sigma}\right)
-\frac{2}{\sqrt{\pi}}e^{-\frac{(t-m)^2}{2\sigma^2}}\sum_{i=0}^{k-1}{k-1
\choose
i}(\sqrt{2}\sigma)^{1-k+i}H_{k-i-1}\left(\frac{x-m}{\sqrt{2}\sigma}\right)\right].
\end{eqnarray*}}On the other hand,$$\frac{d^k}{dt^k}
e^{-\left(\frac{(t-m)^2}{2\sigma^2}+t\right)} =
(-1)^k(t-m+\sigma^2)^k\sigma^{-2k}e^{-\left(\frac{(t-m)^2}{2\sigma^2}+t\right)},$$
and therefore $$e^{-x} \int_{-\infty}^x e^{t}\,\frac{d^k}{dt^k}
e^{-\left(\frac{(t-m)^2}{2\sigma^2}+t\right)}\,\textrm{d}t =
(-1)^k\sigma^{-2k}e^{-x}\int_{-\infty}^{x}(t-m+\sigma^2)^ke^{-\left(\frac{(t-m)^2}{2\sigma^2}+t\right)}\textrm{d}t.$$
Summarizing the above computations we obtain
{\setlength\arraycolsep{2pt}
\begin{eqnarray*}
& {\phantom{1}}&
\int_{-\infty}^{x}(t-m+\sigma^2)^ke^{-\left(\frac{(t-m)^2}{2\sigma^2}+t\right)}\textrm{d}t
\\ & = & \sigma
\sqrt{\frac{\pi}{2}}\left[\textrm{erf}\left(\frac{x-m}{\sqrt{2}\sigma}\right)
-\frac{2}{\sqrt{\pi}}e^{-\frac{(t-m)^2}{2\sigma^2}}\sum_{i=0}^{k-1}{k-1
\choose
i}(\sqrt{2}\sigma)^{1-k+i}H_{k-i-1}\left(\frac{x-m}{\sqrt{2}\sigma}\right)\right].
\end{eqnarray*}}
\end{example}

\begin{example}\rm
Let $\alpha=-\infty$. If $w(t)=te^{-t^2/2}$, then
$w_1(t)/w(t)=-1/t$ for $t>0$. Putting $h(t)=\sin t/t$ we get
$$\left(-\frac{1}{x}\frac{d}{dx}\right)^k\left(e^{\frac{x^2}{2}}\int_{-\infty}^{x} e^{-\frac{t^2}{2}}\sin t\,\textrm{d}t\right)
=
e^{\frac{x^2}{2}}\int_{0}^{x}t^{1-k}e^{-\frac{t^2}{2}}j_k(t)\,\textrm{d}t,$$
where $j_k$ is the spherical Bessel function of the first kind and
the relation
$$\left(-\frac{1}{t}\frac{d}{dt}\right)^k\frac{\sin t}{t} =
\frac{j_k(t)}{t^k}$$ was used (see~\cite{abramovitz}). Since
$$\int_{-\infty}^{x} e^{-\frac{t^2}{2}}\sin t\,\textrm{d}t = \frac{1}{2}\imath\sqrt{\frac{\pi}{2e}}
\left(\textrm{erf}\left(\frac{x+\imath}{\sqrt{2}}\right) -
\textrm{erf}\left(\frac{x-\imath}{\sqrt{2}}\right)\right),$$ where
$\imath$ is the imaginary unit, then we have
$$\int_{0}^{x}t^{1-k}e^{-\frac{t^2}{2}}j_k(t)\,\textrm{d}t = \frac{1}{2}\imath\sqrt{\frac{\pi}{2e}}e^{-\frac{x^2}{2}}
\left(-\frac{1}{x}\frac{d}{dx}\right)^k\left[e^{\frac{x^2}{2}}\left(\textrm{erf}\left(\frac{x+\imath}{\sqrt{2}}\right)
-
\textrm{erf}\left(\frac{x-\imath}{\sqrt{2}}\right)\right)\right].$$
Similarly, if we choose $h(t)=t^{m+1}j_m(t)$, $m$ is integer, and
use the Rayleigh's relation (see~\cite{abramovitz})
$$\left(-\frac{1}{t}\frac{d}{dt}\right)^k(t^{m+1}j_m(t)) =
t^{m-k+1}j_{m-k}(t),$$ the we get
$$\int_{0}^{x}t^{m-k+2}e^{-\frac{t^2}{2}}j_{m-k}(t)\,\textrm{d}t = e^{-\frac{x^2}{2}}
\left(-\frac{1}{x}\frac{d}{dx}\right)^k\left[e^{\frac{x^2}{2}}\int_{0}^{x}t^{m+2}e^{-\frac{t^2}{2}}j_m(t)\,\textrm{d}t\right].$$
\end{example}

\section{Concluding remarks}

As far as authors know many of the above obtained identities does
not appear in the available literature and therefore we suppose
they are new. By all means the method of their acquirement seems
to be new and interesting from the viewpoint of applicability of
some commutation relations for special operators. This method also
provides a way how to obtain other possible identities involving
non-elementary functions (not necessarily special ones) useful in
pure mathematics, theoretical physics as well as applied sciences.



\begin{thebibliography}{99}

\bibitem{abramovitz}
Abramowitz, M., Stegun, I.~A.: Handbook of Mathematical Functions.
Dover, New York (1970).

\bibitem{bell}
Bell, W.~W.: Special Functions for Scientists and Engineers. D.
van Nostrand Comp. Ltd. (1968).

\bibitem{hutnik}
Hutn\'ik, O.: On a general mean-type inequality. In: APLIMAT 2006,
Bratislava, 535--544 (2006).

\bibitem{persson}
Kufner, A., Persson, L.~E.: Weighted Inequalities of Hardy Type.
World Scientific, New Jersey/London/Singapore/Hong Kong (2003).

\bibitem{mojsej}
Mojsej, I.: Asymptotic properties of solutions of third-order
nonlinear differential equations with deviating argument.
Nonlinear Anal. {\bf 68}(11), 3581--3591 (2008).

\bibitem{neumann}
von Neumann, J.: Mathematical Foundations of Quantum Mechanics.
Princeton University Press, New edition (1996).

\bibitem{ohriska}
Ohriska, J.: Oscillation of differential equations and
$v$-derivatives. Czech. Math. J. {\bf 39}(114), 24--44 (1989).

\bibitem{rainville}
Rainville, E.: Special Functions. Macmillan Co., New York (1960),
reprinted by Chelsea Publ. Co., New York (1971).

\end{thebibliography}
\end{document}